\documentclass[11pt]{amsart}
\usepackage{amssymb}
\usepackage{graphicx}
\usepackage[all]{xy}

\title[Flow does not model flows up to weak dihomotopy]{Flow does not model flows up to weak dihomotopy}
\author[P. Gaucher]{Philippe Gaucher}
\address{Preuves Programmes et Syst{\`e}mes\\ Universit{\'e} Paris 7--Denis Diderot\\
Case 7014\\2 Place Jussieu\\ 75251 PARIS Cedex 05\\ France}
\email{gaucher@pps.jussieu.fr}
\urladdr{http://www.pps.jussieu.fr/{\~{}}gaucher/}

\newcommand{\C}{\mathcal{C}}

\newcommand{\p}\times
\renewcommand{\vec}{\overrightarrow}
\renewcommand{\P}{\mathbb{P}}

\newcommand{\be}{\begin{equation}}
\newcommand{\ee}{\end{equation}}
\newcommand{\bea}{\begin{eqnarray}}
\newcommand{\eea}{\end{eqnarray}}
\newcommand{\beas}{\begin{eqnarray*}}
\newcommand{\eeas}{\end{eqnarray*}}

\newtheorem{thm}{Theorem}[section]
\newtheorem*{thmN}{Theorem}
\newtheorem{prop}[thm]{Proposition}
\newtheorem{lem}[thm]{Lemma}

\newtheorem{cor}[thm]{Corollary}
\newtheorem{rem}[thm]{Remark}

\newtheorem{defn}[thm]{Definition}

\newtheorem{nota}[thm]{Notation}
\newcommand{\bd}{\begin{defn}}
\newcommand{\ed}{\end{defn}}
\newcommand{\bcd}{\begin{defn}}
\newcommand{\ecd}{\end{defn}}
\newcommand{\bex}{\begin{exmp}}
\newcommand{\eex}{\end{exmp}}
\newcommand{\bp}{\begin{prop}}
\newcommand{\ep}{\end{prop}}
\newcommand{\bth}{\begin{thm}}
\renewcommand{\eth}{\end{thm}}
\newcommand{\br}{\begin{rem}}
\newcommand{\er}{\end{rem}}
\newcommand{\bpf}{\begin{proof}}
\newcommand{\epf}{\end{proof}}

\newcommand{\fl}[1]{\ar@{->}[l]_{#1}}
\newcommand{\fr}[1]{\ar@{->}[r]^-{#1}}
\newcommand{\fd}[1]{\ar@{->}[d]_{#1}}
\newcommand{\fu}[1]{\ar@{->}[u]^{#1}}
\newcommand{\f}[2]{\ar@{->}[#1]|{#2}}
\newcommand{\ff}[2]{\ar@2{->}[#1]|{#2}}
\newcommand{\frr}[1]{\ar@{->}[rr]^{#1}}

\newcommand{\iso}{\cong}
\newcommand{\lp}{\left(}
\newcommand{\rp}{\right)}

\newcommand{\ot}{\otimes}

\newcommand{\vI}{\vec{I}}

\renewcommand{\geq}{\geqslant}

\def\cartesien{%
  \ar@{-}[]+R+<6pt,-2pt>;[]+RD+<6pt,-6pt>%
  \ar@{-}[]+D+<2pt,-6pt>;[]+RD+<6pt,-6pt>%
}
\def\cocartesien{%
  \ar@{-}[]+L+<-6pt,+2pt>;[]+LU+<-6pt,+6pt>%
  \ar@{-}[]+U+<-2pt,+6pt>;[]+LU+<-6pt,+6pt>%
}

\newcommand{\brm}[1]{\rm{\mathbf{#1}}}

\renewcommand{\top}{{\brm{Top}}}

\newcommand{\dtop}{{\brm{Flow}}}

\newcommand{\set}{{\brm{Set}}}

\newcommand{\glob}{{\rm{Glob}}}

\DeclareMathOperator{\sing}{Sing}

\newcommand{\liminj}{\varinjlim}

\makeatletter

\def\varholim@#1#2{%
  \vtop{\m@th\ialign{##\cr
    \hfil$#1\operator@font holim$\hfil\cr
    \noalign{\nointerlineskip\kern1.5\ex@}#2\cr
    \noalign{\nointerlineskip\kern-\ex@}\cr}}%
}
\def\holimproj{%
  \mathop{\mathpalette\varholim@{\leftarrowfill@\textstyle}}\nmlimits@
}
\def\holiminj{%
  \mathop{\mathpalette\varholim@{\rightarrowfill@\textstyle}}\nmlimits@
}

\makeatother

\DeclareMathOperator{\id}{Id}

\newcommand{\sis}{\Delta^{op}\set}

\DeclareMathOperator{\Iso}{{\underline{Iso}}}
\DeclareMathOperator{\Mono}{{\underline{Mono}}}
\DeclareMathOperator{\Epi}{{\underline{Epi}}}
\DeclareMathOperator{\SplitMono}{{\underline{SplitMono}}}
\DeclareMathOperator{\Empty}{{\underline{Empty}}}
\DeclareMathOperator{\NonEmpty}{{\underline{NonEmpty}}}
\DeclareMathOperator{\All}{{\underline{All}}}
\DeclareMathOperator{\map}{Map}
\DeclareMathOperator{\Mod}{Mod}

\DeclareMathOperator{\cell}{{\brm{cell}}}
\DeclareMathOperator{\cof}{{\brm{cof}}}
\DeclareMathOperator{\inj}{{\brm{inj}}}

\begin{document}

\begin{abstract}
We prove that the category of flows cannot be the underlying category
of a model category whose corresponding homotopy types are the flows
up to weak dihomotopy. Some hints are given to overcome this
problem. In particular, a new approach of dihomotopy involving
simplicial presheaves over an appropriate small category is
proposed. This small category is obtained by taking a full subcategory
of a locally presentable version of the category of flows.
\end{abstract}

\maketitle

\tableofcontents

\section{Introduction}

The category of flows $\dtop$ is introduced in \cite{model3} as a
geometric model of \textit{higher dimensional automata} allowing to
study dihomotopy from the point of view of model category
theory. Roughly speaking, a \textit{weak dihomotopy equivalence} is a
morphism of flows preserving computer scientific properties like the
presence or absence of \textit{deadlocks}, of
\textit{unreachable states} and of \textit{initial and final
states} so that it suffices to work in the categorical
localization. The class of weak dihomotopy equivalences is divided in
two subclasses, the one of \textit{weak S-homotopy equivalences}
\cite{model3} and the one of \textit{T-homotopy equivalences}
\cite{model2}. What is concerned  is the construction of a model
structure (in the sense of \cite{MR36:6480} or \cite{MR99h:55031}) on
the category of flows whose weak equivalences are exactly the weak
dihomotopy equivalences. This way, it becomes possible to study the
categorical localization of the category of flows with respect to weak
dihomotopy equivalence using the tools of algebraic topology.

This is partially done in \cite{model3} where a model structure whose
weak equivalences are exactly the weak S-homotopy equivalences is
constructed. Unfortunately, this model structure does not contain
enough weak equivalences because the T-homotopy equivalences are not
inverted in its homotopy category.

The most elementary example of T-homotopy equivalence which is not
inverted by the model structure constructed in
\cite{model3} is the unique morphism $\phi$ dividing a
directed segment in a composition of two directed segments
(Figure~\ref{ex1} and Notation~\ref{phi})

It is already known that a weak dihomotopy equivalence, whatever it
is, must not be a non-trivial pushout of the morphism of flows
$R:\{0,1\}\longrightarrow \{0\}$ (that is: identifying two distinct
states) because such a pushout either does not preserve one initial or
final state, or creates a loop or a $1$-dimensional branching or a
$1$-dimensional merging (Figure~\ref{loop} and
Figure~\ref{crushing}). The main theorem of this paper is then:

\begin{thmN} (Theorem~\ref{main})
In any model structure on the category of flows such that $\phi$ is a
weak equivalence, there exists a non-trivial pushout of $R$ which is a
weak equivalence. In other terms, there does not exist any model
structure on the category of flows whose weak equivalences are exactly
the weak dihomotopy equivalences.
\end{thmN}

The deep cause of this phenomenon is clearly the presence of
$R:\{0,1\}\longrightarrow \{0\}$ in the class of cofibrations
(Lemma~\ref{main00}). The end of the paper is then devoted to proving
that it is at least possible to get rid of $R$ because:

\begin{thmN} (Theorem~\ref{verscellulaire}) 
Consider the model category of flows whose weak equivalences are
exactly the weak S-homotopy equivalences constructed in
\cite{model3} and recalled in Theorem~\ref{rappel}. 
Then it is Quillen equivalent to another model category whose all
cofibrations are monomorphisms. The new model category contains more
objects but has the same weak S-homotopy types. In particular, the
sets, that is to say the flows with empty path space, are replaced by
the simplicial sets and the epimorphism $R:\{0,1\}\longrightarrow
\{0\}$ by the effective monomorphism $\{0,1\}\subset [0,1]$.
\end{thmN}

Here is now an outline of the paper. Section~\ref{reminder} recalls
what is necessary to know about flows to understand this work.
Section~\ref{restrictionsection} proves that the restriction of any
model structure on the category of flows to the category of sets gives
rise to a model structure on the category of sets. This leads us to
studying in Section~\ref{modelset} the weak factorization systems of
the category of sets.  Then Section~\ref{mainproof} proves the main
theorem of this paper.  At last, Section~\ref{toward} gives some new
directions or research to solve the problem appearing in this paper.

\vspace{0.5cm}
\paragraph{\bf{Acknowledgments}}

I have been trying to construct this model structure on the category
of flows for several months before realizing that it does not
exist. During these months, I had helpful email conversations with
Tibor Beke, Daniel Dugger, Mark Hovey, and especially with Philip
Hirschhorn and I would like to thank them. I also thank Thomas
Goodwillie for his unvolontary help with his remark about the model
structures on the category of sets in Don Davis' mailing-list
``Algebraic Topology''.

\section{Reminder about flows}\label{reminder}

Let $\top$ be the category of compactly generated topological spaces,
i.e.  of weak Hausdorff $k$-spaces. More details for this kind of
topological spaces can be found in \cite{MR90k:54001},
\cite{MR2000h:55002} and the appendix of \cite{Ref_wH}.

\bd\cite{model3}
A {\rm flow} $X$ consists of a topological space $\P X$, a discrete
space $X^0$, two continuous maps $s$ and $t$ called respectively the
source map and the target map from $\P X$ to $X^0$ and a continuous
and associative map $*:\{(x,y)\in \P X\p \P X;
t(x)=s(y)\}\longrightarrow \P X$ such that $s(x*y)=s(x)$ and
$t(x*y)=t(y)$.  A morphism of flows $f:X\longrightarrow Y$ consists of
a set map $f^0:X^0\longrightarrow Y^0$ together with a continuous map
$\P f:\P X\longrightarrow \P Y$ such that $f(s(x))=s(f(x))$,
$f(t(x))=t(f(x))$ and $f(x*y)=f(x)*f(y)$. The corresponding category
is denoted by $\dtop$. \ed

The category $\dtop$ is complete and cocomplete. The topological space
$X^0$ is called the \textit{$0$-skeleton} of $X$. The elements of the
$0$-skeleton $X^0$ are called \textit{states} or \textit{constant
execution paths}.  The elements of $\P X$ are called \textit{non
constant execution paths}. An \textit{initial state} (resp. a
\textit{final state}) is a state which is not the 
target (resp. the source) of any non-constant execution path.

For the sequel, the category of sets $\set$ is identified with the full
subcategory of $\dtop$ consisting of the flows $X$ such that $\P
X=\varnothing$.

\bd\cite{model3}
Let $Z$ be a topological space. Then the {\rm globe} of $Z$ is the
flow $\glob(Z)$ defined as follows: $\glob(Z)^0=\{0,1\}$,
$\P\glob(Z)=Z$, $s=0$, $t=1$ and the composition law is trivial. \ed

\begin{nota} \cite{model3} 
If $Z$ and $T$ are two topological spaces, then $\glob(Z)*\glob(T)$ is
the flow obtained by identifying the final state of $\glob(Z)$ with
the initial state of $\glob(T)$. In other terms, one has the pushout
of flows:
\[
\xymatrix{
\{0\} \fr{0\mapsto 1} \fd{0\mapsto 0} & \glob(Z)\fd{}\\
\glob(T) \fr{} & \cocartesien \glob(Z)*\glob(T)}
\] 
\end{nota}

\begin{nota}\cite{model3} For $\alpha,\beta\in X^0$, 
let $\P_{\alpha,\beta}X$ be the subspace of $\P X$ equipped the
Kelleyfication of the relative topology consisting of the non-constant
execution paths $\gamma$ of $X$ with beginning $s(\gamma)=\alpha$ and
with ending $t(\gamma)=\beta$. \end{nota}

\begin{figure}
\begin{center}
\includegraphics[width=5cm]{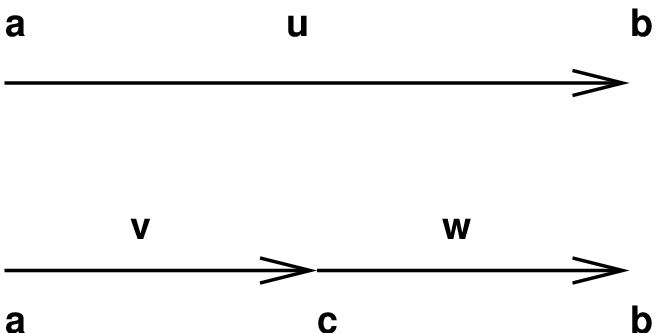}
\end{center}
\caption{Simplest example of T-homotopy equivalence} 
\label{ex1}
\end{figure}

The morphism of flows $\phi$ is going to play an important role in
this paper:

\begin{nota} \label{phi} 
The morphism of flows $\phi:\vI\longrightarrow \vI*\vI$ is the unique
morphism $\phi:\vI\longrightarrow \vI*\vI$ such that
$\phi([0,1])=[0,1]*[0,1]$ where the flow $\vI=\glob(\{[0,1]\})$ is the
directed segment. It corresponds to Figure~\ref{ex1}. \end{nota}

\begin{figure}
\begin{center}
\includegraphics[width=5cm]{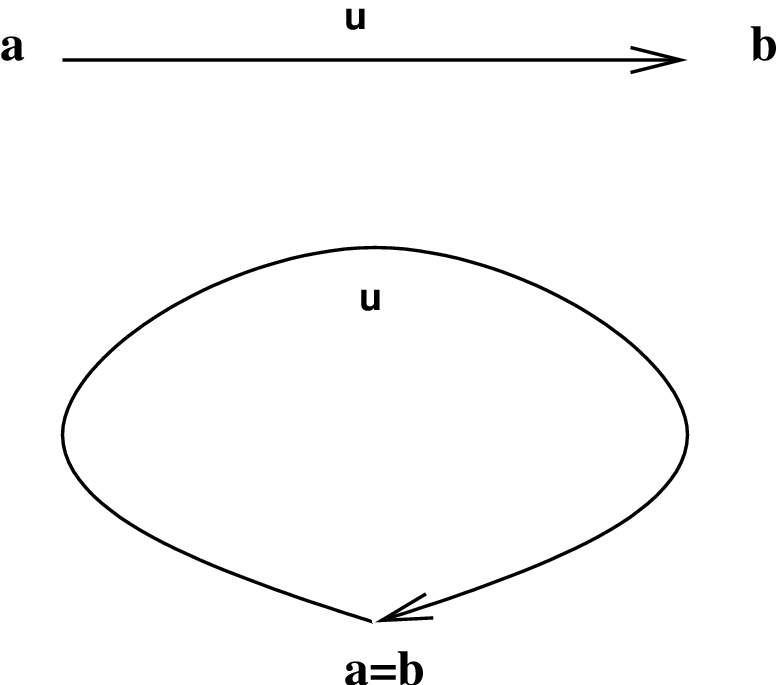}
\end{center}
\caption{Non-authorized identification} 
\label{loop}
\end{figure}

\begin{figure}
\begin{center}
\includegraphics[width=5cm]{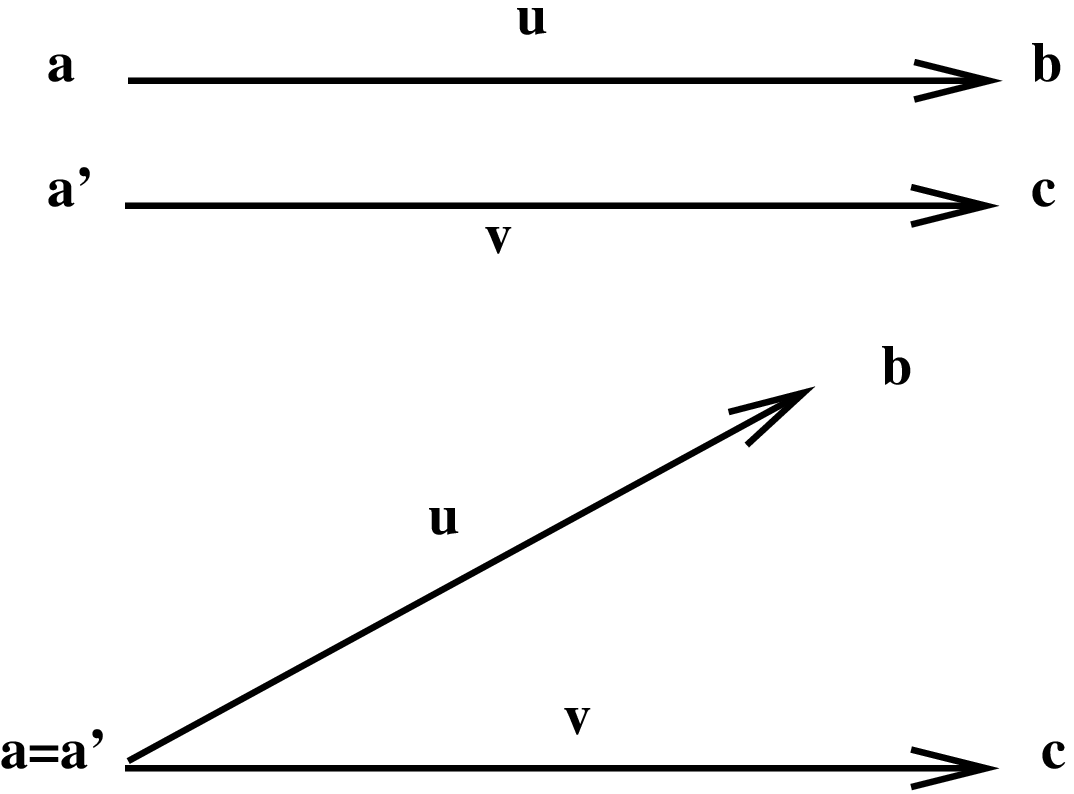}
\end{center}
\caption{Non-authorized identification} 
\label{crushing}
\end{figure}

The morphism of flows $\phi:\vI\longrightarrow \vI*\vI$ is an example
of T-homotopy equivalence, as introduced in \cite{model2}. We would
want the morphism $\phi$ to be a weak equivalence since it is an
example of refinement of observation. On the contrary, the non-trivial
pushouts of the morphism $R:\{0,1\}\longrightarrow
\{0\}$ must absolutely be removed from the class of weak equivalences because any
such pushout either does not preserve one initial or final state, or
creates a loop or a $1$-dimensional branching or a $1$-dimensional
merging (Figure~\ref{loop} and Figure~\ref{crushing}). So such a
pushout cannot be a dihomotopy equivalence.

\section{Restriction to $\set$ of a model structure on $\dtop$}
\label{restrictionsection}

For any category $\C$, $\map(\C)$ denotes the class of morphisms of
$\C$. In a category $\C$, an object $x$ is \textit{a retract} of an
object $y$ if there exist $f:x\longrightarrow y$ and
$g:y\longrightarrow x$ of $\C$ such that $g\circ f=\id_x$. A
\textit{functorial factorization} $(\alpha,\beta)$ of $\C$ is a pair
of functors from $\map(\C)$ to $\map(\C)$ such that for any $f$ object
of $\map(\C)$, $f=\beta(f)\circ \alpha(f)$.

\bd Let $i:A\longrightarrow B$ and $p:X\longrightarrow Y$ be maps in a
category $\C$. Then $i$ has the {\rm left lifting property} (LLP)
with respect to $p$ (or $p$ has the {\rm right lifting property}
(RLP) with respect to $i$) if for any commutative square
\[
\xymatrix{
A\fd{i} \fr{\alpha} & X \fd{p} \\
B \ar@{-->}[ru]^{g}\fr{\beta} & Y}
\]
there exists $g$ making both triangles commutative. \ed

\bd \cite{MR2003h:18001}
Let $\C$ be a category. A {\rm weak factorization system} is a pair
$(\mathcal{L},\mathcal{R})$ of classes of morphisms of $\C$ such that
the class $\mathcal{L}$ is the class of morphisms having the LLP with
respect to $\mathcal{R}$, such that the class $\mathcal{R}$ is the
class of morphisms having the RLP with respect to $\mathcal{L}$ and
such that any morphism of $\C$ factors as a composite $r\circ \ell$
with $\ell\in \mathcal{L}$ and $r\in \mathcal{R}$. The weak
factorization system is {\rm functorial} if the factorization $r\circ
\ell$ is a functorial factorization. \ed

In a weak factorization system $(\mathcal{L},\mathcal{R})$, the class
$\mathcal{L}$ (resp.  $\mathcal{R}$) is completely determined by
$\mathcal{R}$ (resp.  $\mathcal{L}$).

\bd \cite{MR99h:55031} 
A {\rm model category} is a complete and cocomplete category equipped
with three classes of morphisms $({\rm Cof},{\rm Fib},\mathcal{W})$
(resp. called the classes of cofibrations, fibrations and weak
equivalences) such that:
\begin{enumerate}
\item the class of morphisms $\mathcal{W}$ is closed under retracts and 
satisfies the two-out-of-three
axiom i.e.:  
if $f$ and $g$ are morphisms of $\C$ such that $g\circ f$
is defined and two of $f$, $g$ and $g\circ f$ are weak
equivalences, then so is the third.
\item the pairs $({\rm Cof}\cap \mathcal{W}, {\rm Fib})$ and
$({\rm Cof}, {\rm Fib}\cap \mathcal{W})$ are both functorial weak factorization
systems.
\end{enumerate}
The triple $({\rm Cof},{\rm Fib},\mathcal{W})$ is called a {\rm
model structure}. An element of ${\rm Cof}\cap \mathcal{W}$ is called 
a {\rm trivial cofibration}. An element of ${\rm Fib}\cap \mathcal{W}$ is called 
a {\rm trivial fibration}.
\ed

\begin{lem}\label{nonvide}
If $f:X\longrightarrow Y$ is a morphism of flows such that either the
space $\P X$ or the space $\P Y$ is non-empty, then $f$ satisfies the
LLP with respect to any set map.
\end{lem}

\bpf This is due to the fact that there does not exist any continuous map
from a non-empty space to an empty space.
\epf

\bth\label{restriction}
Let $({\rm Cof},{\rm Fib},\mathcal{W})$ be a model structure on
$\dtop$.  Then \[({\rm Cof}\cap \map(\set),{\rm Fib}\cap
\map(\set),\mathcal{W}\cap \map(\set))\] is a model structure on
$\set$. \eth

\bpf Let $f:X\longrightarrow Y$ and $g:Y\longrightarrow Z$ be two morphisms of
flows. If $f$ and $g$ belong to $\mathcal{W}\cap \map(\set)$, then $\P
X=\P Y=\P Z=\varnothing$, and so $g\circ f\in \mathcal{W}\cap
\map(\set)$. If $f$ and $g\circ f$ belong to $\mathcal{W}\cap
\map(\set)$, then $\P X=\P Y=\varnothing$ and $\P Z=\varnothing$. So
$g\in \mathcal{W}\cap \map(\set)$. If $g$ and $g\circ f$ belong to
$\mathcal{W}\cap \map(\set)$, then $\P Y=\P Z=\varnothing$ and $\P
X=\varnothing$. So $f\in \mathcal{W}\cap \map(\set)$.  The class
$\mathcal{W}\cap \map(\set)$ is closed under retracts since the only
retract of the empty set is the empty set. So $\mathcal{W}\cap
\map(\set)$ satisfies the two-out-of-three axiom and is closed under
retracts.

Let $g\in \map(\set)$ satisfying the RLP with respect to any morphism
of ${\rm Cof}\cap \map(\set)$. By Lemma~\ref{nonvide}, $g$ satisfies
the RLP with respect to any morphism of $\map(\dtop)\backslash
\map(\set)$.  But ${\rm Cof}\subset \lp\map(\dtop)\backslash
\map(\set)\rp\cup \lp {\rm Cof}\cap \map(\set)\rp$. So
$g\in {\rm Fib}\cap\mathcal{W}\cap\map(\set)$. Thus, the class ${\rm
Fib}\cap\mathcal{W}\cap\map(\set)$ is exactly the class of morphisms
of sets satisfying the RLP with respect to ${\rm Cof}\cap
\map(\set)$. Let $h\in \map(\set)$. Then $h=r\circ i$ with
$i\in {\rm Cof}$ and $r\in {\rm Fib}\cap\mathcal{W}$. So $r$ induces
a continuous map from the path space of its domain to the empty space.
So both $r$ and $i$ are set maps. Therefore the pair
$({\rm Cof}\cap\map(\set), {\rm Fib}\cap\mathcal{W}\cap\map(\set))$
is a functorial weak factorization system.

And for similar reasons, the pair $({\rm
Cof}\cap\mathcal{W}\cap\map(\set), {\rm Fib}\cap\map(\set))$ is a
functorial weak factorization system as well. \epf

\section{The weak factorization systems of the category of sets}
\label{modelset}

Theorem~\ref{restriction} leads us to studying the possible weak
factorization systems of the category of sets. Let us describe them
now.

Let $\All$ be the class of all morphisms of sets. Let $\Iso$ be the
class of bijections. Let $\Mono$ be the class of injections.  Let
$\Epi$ be the class of surjections. Let $\SplitMono$ be the class of
set maps having a left inverse. Let $\Empty$ be the class of set maps
whose domain is the empty set. Let $\NonEmpty$ be the class of set
maps whose domain is non-empty.

Let $C:\varnothing\longrightarrow \{0\}$. Let
$C^+:\{0\}\longrightarrow \{0,1\}$ with $C^+(0)=0$. Let
$R:\{0,1\}\longrightarrow \{0\}$. One has $C\in\Empty\subset\Mono$,
$C\notin \SplitMono$, $C^+\in\SplitMono\subset\Mono$ and $R\in\Epi$.

\begin{nota} Let $\C$ be a cocomplete category. 
If $K$ is a set of morphisms of $\C$, then the collection of
morphisms of $\C$ that satisfy the RLP with respect to any morphism
of $K$ is denoted by $\inj(K)$ and the collection of morphisms of $\C$
that are transfinite compositions of pushouts of elements of $K$ is
denoted by $\cell(K)$. Denote by $\cof(K)$ the collection of morphisms
of $\C$ that satisfies the LLP with respect to any morphism that
satisfies the RLP with respect to any element of $K$.  This is a
purely categorical fact that $\cell(K)\subset \cof(K)$.  \end{nota}

\begin{lem} For any set of morphisms of sets $K$, the pair $(\cof(K),\inj(K))$
is a functorial weak factorization system. \end{lem}

\bpf This is due to the fact that any set is small (in the sense of model categories) 
so the small object argument applies \cite{MR99h:55031}. This is also
a consequence of the fact that $\set$ is locally presentable
\cite{MR95j:18001} and of \cite{MR1780498} Proposition~1.3. \epf

\begin{lem}\label{cons} 
For any set of morphisms of sets $K$, the class of morphisms of sets
$\cof(K)$ is exactly the class of retracts of transfinite compositions
of pushouts of elements of $K$. \end{lem}

\bpf This is due to the fact that $\set$ is cocomplete, that 
any set is small and to  \cite{MR99h:55031} Corollary~2.1.15. \epf

\begin{lem} \label{allwfs0}
One has:
\begin{enumerate}
\item if $K=\varnothing$, then $(\cof(K),\inj(K))=(\Iso,\All)$
\item if $K=\{C\}$, then $(\cof(K),\inj(K))=(\Mono,\Epi)$
\item if $K=\{C^+\}$, then $(\cof(K),\inj(K))=(\SplitMono, \Epi\cup \Empty)$
\item if $K=\{R\}$, then $(\cof(K),\inj(K))=(\Epi,\Mono)$
\item if $K=\{R,C\}$, then $(\cof(K),\inj(K))=(\All,\Iso)$
\item if $K=\{R,C^+\}$, then $(\cof(K),\inj(K))=(\Iso \cup \NonEmpty , \Iso
\cup \Empty)$.
\end{enumerate}
\end{lem}

\bpf By Lemma~\ref{cons}, the class $\cof(K)$ is exactly the class of 
retracts of transfinite compositions of pushouts of maps of $K$. Hence
$\cof(\varnothing)=\Iso$, $\cof(\{C\})=\Mono$, $\cof(\{C^+\})=\SplitMono$,
$\cof(\{R\})=\Epi$, $\cof(\{R,C\})=\All$ and $\cof(\{R,C^+\})=\Iso \cup
\NonEmpty$. The equalities $\inj(\varnothing)=\All$, 
$\inj(\{C\})=\Epi$ and $\inj(\{R\})=\Mono $ are clear. The equality
$\inj(\{C^+\})=\Epi\cup \Empty$ is a consequence of $\inj(\{C\})=\Epi$ and
of the fact that there does not exist any set map from a non-empty set
to the empty set. And $\inj(\{R,C\})=\inj(\{R\})\cap
\inj(\{C\})=\Mono\cap\Epi=\Iso$.  At last: $\inj(\{R,C^+\})=\inj(\{R\})\cap
\inj(\{C^+\})=\Mono\cap (\Epi\cup \Empty)=\Iso
\cup \Empty$. \epf

\bth\label{allwfs} (Goodwillie)
The six weak factorization systems $(\Iso,\All)$, $(\Mono,\Epi)$,
$(\SplitMono, \Epi\cup \Empty)$, $(\Epi,\Mono)$, $(\All,\Iso)$ and
$(\Iso \cup \NonEmpty , \Iso\cup \Empty)$ are the only possible weak
factorization systems on the category of sets.
\eth

\bpf Let $(\mathcal{L},\mathcal{R})$ be a weak factorization system on the
category of sets. Then $\Iso\subset \mathcal{L}\cap \mathcal{R}$.

First of all, if $\mathcal{L}=\Iso$, then $\mathcal{R}=\All$ by
Lemma~\ref{allwfs0}. Let us suppose now that $\mathcal{L}\backslash
\Iso\neq \varnothing$.

If $C\in \mathcal{L}$, then $\Mono\subset \mathcal{L}$ since
$\mathcal{L}$ is closed under pushouts and transfinite composition. So
if $\mathcal{L}\subset \Mono$, then $\mathcal{L}=\Mono$ and
necessarily $(\mathcal{L},\mathcal{R})=(\Mono,\Epi)$ by
Lemma~\ref{allwfs0}.

If $C\in \mathcal{L}$ and if $\mathcal{L}{\not\subset} \Mono$, let
$f\in \mathcal{L}\backslash \Mono$. Then $R$ is a retract of $f$ and
so $R\in\mathcal{L}$. Therefore in this case, $\mathcal{L}=\All$ since
$\mathcal{L}$ is closed under pushouts and transfinite compositions
and necessarily $\mathcal{R}=\Iso$ by Lemma~\ref{allwfs0}.

Let us suppose now that $C\notin \mathcal{L}$ and that $C^+\in
\mathcal{L}$.  Then $\SplitMono \subset \mathcal{L}$ since
$\mathcal{L}$ is closed under pushouts and transfinite composition.
So if $\mathcal{L}\subset \Mono$, then $\mathcal{L}=\SplitMono$ and
$\mathcal{R}=\Epi\cup \Empty$ by Lemma~\ref{allwfs0}. And if
$\mathcal{L}{\not\subset}\Mono$, then $R\in \mathcal{L}$ like
above. So $\Iso \cup \NonEmpty\subset \mathcal{L}$ since $\mathcal{L}$
is closed under pushouts and transfinite composition. And therefore
since $C\notin \mathcal{L}$, one has $\Iso \cup \NonEmpty=
\mathcal{L}$ and $\mathcal{R}=\Iso
\cup \Empty$ by Lemma~\ref{allwfs0}.

Let us suppose now that $C\notin \mathcal{L}$ and that $C^+\notin
\mathcal{L}$. Then $R\in \mathcal{L}$ like above. And
$\Epi\subset \mathcal{L}$ since $\mathcal{L}$ is closed under pushouts
and transfinite composition. So in this case, one has
$\mathcal{L}=\Epi$ and $\mathcal{R}=\Mono$ by Lemma~\ref{allwfs0}. \epf

It is even possible to prove the:

\bth (Goodwillie's exercise) The nine model structures of the category of sets are: 
\begin{enumerate}
\item $({\rm Cof},{\rm Fib},\mathcal{W})=(\All ,\All ,\Iso )$
\item $({\rm Cof},{\rm Fib},\mathcal{W})=(\All ,\Iso \cup\Empty , \Iso \cup\NonEmpty)$
\item $({\rm Cof},{\rm Fib},\mathcal{W})=( \All, \Iso, \All)$
\item $({\rm Cof},{\rm Fib},\mathcal{W})=( \Iso, \All, \All)$
\item $({\rm Cof},{\rm Fib},\mathcal{W})=( \Epi, \Mono,\All )$
\item $({\rm Cof},{\rm Fib},\mathcal{W})=( \Mono,\Epi , \All)$
\item $({\rm Cof},{\rm Fib},\mathcal{W})=( \SplitMono,\Epi\cup \Empty ,\All )$
\item $({\rm Cof},{\rm Fib},\mathcal{W})=(\Iso\cup \NonEmpty , \Iso\cup\Empty, \All)$
\item $({\rm Cof},{\rm Fib},\mathcal{W})=(\Mono , \Epi\cup \Empty,\Iso \cup\NonEmpty )$. 
\end{enumerate}
\eth

\bpf Exercise for the idle mathematician proposed by Goodwillie in 
Don Davis' mailing-list ``Algebraic Topology''. \epf

Only Theorem~\ref{allwfs} is necessary for the proof of the main
theorem.

\section{Proof of the main theorem}\label{mainproof}

\begin{lem} \label{c1}
Let $({\rm Cof},{\rm Fib},\mathcal{W})$ be a model structure on
$\dtop$ whose weak equivalences are never a non-trivial pushout of
$R$. Then the only possibilities for the weak factorization system
$({\rm Cof}\cap
\mathcal{W}\cap\map(\set),{\rm Fib}\cap\map(\set))$ are
$(\Iso,\All)$, $(\Mono,\Epi)$ and $(\SplitMono, \Epi\cup \Empty)$.
\end{lem}

\bpf The morphism $R$ is not a weak equivalence, so it cannot 
be a trivial cofibration. So $R\notin {\rm Cof}\cap\mathcal{W}\cap
\map(\set)$. The proof is then complete with
Theorem~\ref{allwfs}. \epf

\begin{lem} \label{c2}
Let $({\rm Cof},{\rm Fib},\mathcal{W})$ be a model structure on
$\dtop$ whose weak equivalences are never a non-trivial pushout of
$R$.  Then the only possibilities for the weak factorization system
$({\rm Cof}\cap\map(\set),{\rm Fib}\cap
\mathcal{W}\cap\map(\set))$ are $(\Epi,\Mono)$, $(\All,\Iso)$ and
$(\Iso \cup\NonEmpty , \Iso\cup \Empty)$.
\end{lem}

\bpf The morphism $R$ is not a weak equivalence, so it cannot 
be a trivial fibration. So $R\notin {\rm Fib}\cap\mathcal{W}\cap
\map(\set)$.  The proof is then complete with
Theorem~\ref{allwfs}. \epf

\begin{lem} \label{main00}
Let $({\rm Cof},{\rm Fib},\mathcal{W})$ be a model structure on
$\dtop$ whose weak equivalences are never a non-trivial pushout of
$R$.  Then both $C:\varnothing\longrightarrow \{0\}$ and
$R:\{0,1\}\longrightarrow \{0\}$ are cofibrations. \end{lem}

\bpf By Theorem~\ref{restriction}, the triple 
$({\rm Cof}\cap\map(\set),{\rm
Fib}\cap\map(\set),\mathcal{W}\cap\map(\set))$ yields a model
structure on the category of sets.  By Lemma~\ref{c1} and
Lemma~\ref{c2}, we then have $3\p 3=9$ possibilities for this
restriction. These nine possibilities are summarized in
Table~\ref{finalement}.

\begin{table}
\begin{center}
{\small
\begin{tabular}{|c|c|c|c|}
\hline
& $(\Epi,\Mono)$ & $(\All,\Iso)$ & $(\Iso \cup\NonEmpty , \Iso\cup \Empty)$ \\
\hline
$(\Iso,\All)$ & $\mathcal{W}=\Mono$ & possible &$\mathcal{W}=\Iso\cup \Empty$\\
\hline
$(\Mono,\Epi)$ & $\Mono\not\subset \Epi$  &$\mathcal{W}=\Mono$& 
$\Iso\cup \Empty\not\subset\Epi$ \\
\hline
$(\SplitMono, \Epi\cup \Empty)$ & $\SplitMono\not\subset\Epi$ &
$\mathcal{W}=\SplitMono$& $\Iso\cup
\Empty\subset \mathcal{W}$\\
\hline
\end{tabular}
}
\end{center}
\caption{The last nine   possibilities}
\label{finalement}
\end{table}

Three situations are impossible because the class of trivial
cofibrations (resp. of trivial fibrations) must be a subclass of the
one of cofibrations (resp. of fibrations): $\Mono\not\subset \Epi$,
$\Iso\cup \Empty\not\subset\Epi$ and $\SplitMono\not\subset\Epi$.

Four other situations are impossible since the class $\mathcal{W}$ of
weak equivalences cannot satisfy the two-out-of-three axiom:
$\mathcal{W}=\Mono$ (twice), $\mathcal{W}=\Iso\cup \Empty$ and
$\mathcal{W}=\SplitMono$.

The column $(\Iso \cup\NonEmpty , \Iso\cup \Empty)$ implies that the
class $\mathcal{W}$ of weak equivalences satisfies $\Iso\cup
\Empty\subset \mathcal{W}$. Consider the composite
$\varnothing\longrightarrow X\longrightarrow Y$ for any set map from
$X$ to $Y$. One deduces that $\mathcal{W}=\All$.  So
$\SplitMono=\All\cap(\Iso \cup\NonEmpty)$: contradiction.

Therefore it remains the case $(\Iso,\All),(\All,\Iso)$ which does
correspond to a possible model structure for the category of sets,
that is $({\rm Cof},{\rm Fib},\mathcal{W})=(\All,\All,\Iso)$.
\epf

\begin{lem}\label{main0}
Let $({\rm Cof},{\rm Fib},\mathcal{W})$ be a model structure on
$\dtop$ whose weak equivalences are never a non-trivial pushout of
$R$. Then any trivial fibration induces a bijection between the
$0$-skeletons.
\end{lem}

\bpf By Lemma~\ref{main00}, any trivial fibration $r$
satisfies the RLP with respect to $R$ and $C$. Thus, the set map $r^0$
is bijective. \epf

\begin{lem}\label{main1} 
Let $({\rm Cof},{\rm Fib},\mathcal{W})$ be a model structure on
$\dtop$ whose weak equivalences are never a non-trivial pushout of
$R$. Then for any trivial cofibration $f$, the set map $f^0$ is
one-to-one. \end{lem}

\bpf Let us suppose that there exists a trivial cofibration
$f:X\longrightarrow Y$ such that $f^0$ is not one-to-one, that is
there exists $(\alpha,\beta)\in X^0\p X^0$ such that $\alpha\neq \beta$
and $f^0(\alpha)=f^0(\beta)$. Then consider the diagram of flows
\[
\xymatrix{
\{0,1\} \fd{\iota} \fr{R} & \{0\}\fd{}\\
X\fd{f} \fr{\widehat{R}} & Z\cocartesien \fd{\widehat{f}}\\
Y \ar@{=}[r] & Y\cocartesien}
\]
with $\iota(0)=\alpha$, $\iota(1)=\beta$. Since $f$ is a trivial
cofibration, $\widehat{f}$ is a trivial cofibration as well. Since
$\widehat{f}\circ \widehat{R}$ is a weak equivalence, the morphism
$\widehat{R}$ therefore belongs to $\mathcal{W}$. Contradiction. \epf

\bth
Let $({\rm Cof},{\rm Fib},\mathcal{W})$ be a model structure on
$\dtop$ whose weak equivalences are never a non-trivial pushout of
$R$. Then for any $f\in \mathcal{W}$, $f^0$ is a bijection.
\eth

\bpf Let $g:X\longrightarrow Y$ of ${\rm Cof}\cap \mathcal{W}$.
By Lemma~\ref{main1}, $g^0$ is a one-to-one set map. Let us
suppose that $g^0$ is not bijective. Consider the diagram of flows
\[
\xymatrix{
X\fr{g} \fd{g} & Y \fd{k}\ar@/^15pt/[rdd]^{\id_Y}&\\ Y
\ar@/_15pt/[rrd]_{\id_Y}\fr{k} &\cocartesien Y\sqcup_X
Y\ar@{-->}[rd]^{h} &\\ &&Y}
\]
Since $g$ is a trivial cofibration of $({\rm Cof},{\rm
Fib},\mathcal{W})$, the morphism of flows $k$ is a trivial cofibration
of $({\rm Cof},{\rm Fib},\mathcal{W})$ as well. Since $h\circ k$ is an
isomorphism and therefore a weak equivalence, $h$ is therefore a weak
equivalence of $({\rm Cof},{\rm Fib},\mathcal{W})$. Moreover, $h^0$ is
epi and not bijective since $g^0$ is not bijective and since one has
the diagram of sets
\[
\xymatrix{
X^0\fr{g^0} \fd{g^0} & Y^0 \fd{k^0}\ar@/^15pt/[rdd]^{\id_{Y^0}}&\\ Y^0
\ar@/_15pt/[rrd]_{\id_{Y^0}}\fr{k^0} &\cocartesien Y^0\sqcup_{X^0}
Y^0\ar@{-->}[rd]^{h^0} &\\ &&Y^0}
\]
The morphism of
flows $h$ factors as $h=r\circ i$ with $i\in {\rm Cof}\cap
\mathcal{W}$ and $r\in {\rm Fib}$. By Lemma~\ref{main1} again, $i^0$
is one-to-one. Since $h=r\circ i\in
\mathcal{W}$, the morphism of flows $r$ belongs to $\mathcal{W}$ as
well. So $r\in {\rm Fib}\cap
\mathcal{W}$. By Lemma~\ref{main0}, $r^0$ is
bijective. Therefore $h^0$ is a one-to-one set map. Contradiction.  So
$g^0$ is bijective.

Any morphism $f$ of $\mathcal{W}$ factors as a composite $f=q\circ j$
with $q\in {\rm Fib}\cap\mathcal{W}$ and $j\in
{\rm Cof}\cap\mathcal{W}$. Since $q^0$ and $j^0$ are both bijective,
$f^0=q^0\circ j^0$ is bijective as well.
\epf

The weak S-homotopy model structure of flows constructed in
\cite{model3} and the Cole-Str{\o}m model structure of flows
constructed in \cite{math.AT/0401033} are two examples of model
category structure on $\dtop$ such that any weak equivalence is never
a non-trivial pushout of $R$.

\bth  \label{main}
In any model structure on the category of flows containing
$\phi:\vI\longrightarrow \vI*\vI$ as weak equivalence, there exists a
non-trivial pushout of $R$ which is a weak equivalence. In other
terms, there does not exist any model structure on the category of
flows whose weak equivalences are exactly the weak dihomotopy
equivalences.
\eth

\bpf The morphism of flows $\phi:\vI\longrightarrow \vI*\vI$ does
not induce a bijection between $\vI^0=\{0,1\}$ and the $0$-skeleton
of $\vI*\vI$ which is a $3$-element set. \epf

\begin{cor}
For any model structure on the category of flows such that $\phi$ is a
weak equivalence, there exists a weak equivalence which does not
preserve the branching homology or the merging homology.
\end{cor}

\section{Towards other models for dihomotopy}
\label{toward}

Theorem~\ref{main} shows that the category of flows cannot be the
underlying category of a model category whose corresponding homotopy
types are the flows up to weak dihomotopy. The cause of the problem
seems to be the unavoidable presence of $R:\{0,1\}\longrightarrow
\{0\}$ in the class of cofibrations (Proposition~\ref{main00}). 
In particular, it prevents the weak S-homotopy model structure of
$\dtop$ from being cellular in the sense of \cite{ref_model2}.

Let $n\geq 1$. Let $\mathbf{D}^n$ be the closed $n$-dimensional disk
and let $\mathbf{S}^{n-1}$ be its boundary. Let
$\mathbf{D}^{0}=\{0\}$. Let $\mathbf{S}^{-1}=\varnothing$ be the empty
space. Let us recall the:

\bth \cite{model3} \label{rappel}
The category of flows $\dtop$ is given a structure of cofibrantly 
generated model category such that: 
\begin{enumerate}
\item the set of generating cofibrations is the union of $\{R,C\}$ and 
the set of morphisms $\glob(f)$ for $f$ running over the set of inclusions
$\glob(\mathbf{S}^{n-1})\longrightarrow \glob(\mathbf{D}^{n})$ for
$n\geq 0$
\item the set of generating trivial cofibrations is the set 
of morphisms $\glob(f)$ for $f$ running over the set of inclusions
$\glob(\mathbf{D}^{n}\p\{0\})\longrightarrow
\glob(\mathbf{D}^{n}\p[0,1])$ 
\item a morphism $f:X\longrightarrow Y$ of $\dtop$ is a weak equivalence 
if and only if $f:X^0\longrightarrow Y^0$ is a bijection of sets and
for any $(\alpha,\beta)\in X^0\p X^0$,
$f:\P_{\alpha,\beta}X\longrightarrow \P_{f(\alpha),f(\beta)}Y$ is a
weak homotopy equivalence of topological spaces, that is a {\rm weak
S-homotopy equivalence}.
\end{enumerate} 
This model structure is called the {\rm weak S-homotopy model structure} of
$\dtop$.  In this model structure, any object is fibrant.
\eth

We are going to prove in this section the:

\bth \label{verscellulaire} 
The model category $\dtop$ is Quillen equivalent to a model category
whose all cofibrations are monomorphisms (and even effective
monomorphisms in the sense of \cite{ref_model2}, or regular
monomorphisms in the sense of \cite{MR96g:18001a} \cite{MR96g:18001b}).
\eth

For this purpose, let us introduce the notion of \textit{flow over a monoidal
category}. Let $(\C,\ot)$ be a monoidal category.

An object $X$ of
$\dtop(\C,\ot)$ consists of a set $X^0$ called the
\textit{$0$-skeleton} of $X$ and for any $(\alpha,\beta)\in X^0\p X^0$
an object $\P_{\alpha,\beta}X$ of $\C$ such that there exists a
morphism $*:\P_{\alpha,\beta}X\ot
\P_{\beta,\gamma}X\longrightarrow \P_{\alpha,\gamma}X$ of $\C$ for any
$(\alpha,\beta,\gamma)\in X^0\p X^0\p X^0$ satisfying the
associativity axiom: for any $(\alpha,\beta,\gamma,\delta)\in X^0\p
X^0\p X^0\p X^0$, the following diagram is commutative
\[
\xymatrix{
\P_{\alpha,\beta}X\ot \P_{\beta,\gamma}X\ot \P_{\gamma,\delta}X \fr{(*,\id)}\fd{(\id,*)} & \P_{\alpha,\gamma}X\ot \P_{\gamma,\delta}X \fd{*} \\
\P_{\alpha,\beta}X\ot \P_{\beta,\delta}X \fr{*}& \P_{\alpha,\delta}X}
\]

A morphism $f:X\longrightarrow Y$ of $\dtop(\C,\ot)$ consists of a set
map $f:X^0\longrightarrow Y^0$ together with morphisms
$f:\P_{\alpha,\beta}X\longrightarrow \P_{f(\alpha),f(\beta)}Y$ such
that the following diagram is commutative for any
$(\alpha,\beta,\gamma)\in X^0\p X^0\p X^0$
\[
\xymatrix{
\P_{\alpha,\beta}X\ot \P_{\beta,\gamma}X\fr{(f,f)} \fd{*}& \P_{f(\alpha),f(\beta)}Y\ot \P_{f(\beta),f(\gamma)}Y \fd{*}\\
\P_{\alpha,\gamma}X \fr{f} & \P_{f(\alpha),f(\gamma)}Y}
\]

\begin{nota} Let $Z$ be an object of $\C$. Denote by $\glob(Z)$ the flow such
that $\glob(Z)^0=\{0,1\}$ and $\P_{0,1}\glob(Z)=Z$. This defines a
full and faithful functor $\glob:\C\longrightarrow
\dtop(\C,\ot)$. \end{nota}

\begin{nota} From now on, the category $\dtop$ is denoted by $\dtop(\top,\p)$. 
\end{nota}

\begin{nota} 
The pair $(\sis,\p)$ denotes the monoidal model category of
simplicial sets \cite{MR2001d:55012}. \end{nota}

\begin{nota} The geometric realization functor is denoted by 
$|-|:\sis\longrightarrow \top$. The singular nerve functor is denoted
by $\sing:\top\longrightarrow \sis$.
\end{nota}

For the sequel, the categories $\sis$ and $\top$ are supposed to be
equipped with their standard cofibrantly generated model structure.

The following two lemmas are necessary for the proof of
Theorem~\ref{modelsis}.

\begin{lem} \label{llp1}
Let $f:U\longrightarrow V$ be a morphism of simplicial sets. Then
\[\glob(f):\glob(U)\longrightarrow
\glob(V)\] satisfies the LLP with respect to the morphism
$g:X\longrightarrow Y$ of $\dtop(\sis,\p)$ if and only if for any
$(\alpha,\beta)\in X^0\p X^0$, $f$ satisfies the LLP with respect to
$\P_{\alpha,\beta}X\longrightarrow\P_{g(\alpha),g(\beta)}Y$. \end{lem}

\bpf Obvious. \epf

The category of sets can be viewed as a full subcategory of the
category of flows over simplicial sets by identifying a set $X$ with
the flow $Y$ such that $Y=X$ and $\P Y=\varnothing$.

\begin{lem} \label{llp2} 
Let $f:U\longrightarrow V$ be a set map. Then $f:U\longrightarrow V$
satisfies the LLP with respect to the morphism $g:X\longrightarrow Y$
of $\dtop(\sis,\p)$ if and only if $f$ satisfies the LLP with respect
to $g^0$. \end{lem}

\bpf Obvious. \epf

\bth 
\label{modelsis} 
The category of flows $\dtop(\sis,\p)$ is given a structure of cofibrantly 
generated model category such that:
\begin{enumerate}
\item the set $\mathcal{I}$ of generating cofibrations is the union of 
$\{R,C\}$ and the set of morphisms $\glob(f)$ for $f$ running over the
set of generating cofibrations of the cofibrantly generated model
category $\sis$ of simplicial sets
\item the set $\mathcal{J}$ of generating trivial cofibrations is the set 
of morphisms $\glob(f)$ for $f$ running over the set of generating
trivial cofibrations of the cofibrantly generated model category
$\sis$ of simplicial sets
\item a morphism $f:X\longrightarrow Y$ of $\dtop(\sis,\p)$ is a 
weak equivalence if and only if $f:X^0\longrightarrow Y^0$ is a
bijection of sets and for any $(\alpha,\beta)\in X^0\p X^0$,
$f:\P_{\alpha,\beta}X\longrightarrow \P_{f(\alpha),f(\beta)}Y$ is a
weak homotopy equivalence of simplicial sets.
\end{enumerate} 
Moreover, the weak S-homotopy model category $\dtop(\top,\p)$ of
\cite{model3} is Quillen equivalent to the model category
$\dtop(\sis,\p)$.
\eth

\bpf 
The class of weak equivalences clearly satisfies the two-out-of-three
axiom.  Any object of $\dtop(\sis,\p)$ is small since any simplicial
set is small by \cite{MR99h:55031} Lemma~3.1.1. So we only have to
check that $\cell(\mathcal{J})\subset \cof(\mathcal{I})\cap
\mathcal{W}$ where $\mathcal{W}$ denotes the class of weak
equivalences and that
$\inj(\mathcal{I})=\inj(\mathcal{J})\cap\mathcal{W}$ by
\cite{MR99h:55031} Theorem~2.1.19.

Let $g:X\longrightarrow Y\in
\inj(\mathcal{I})$. Then for any $(\alpha,\beta)\in X^0\p X^0$, the morphism 
of simplicial sets $\P_{\alpha,\beta}X\longrightarrow
\P_{g(\alpha),g(\beta)}Y$ satisfies the RLP with respect to any
cofibration of simplicial sets by Lemma~\ref{llp1}. So
$\P_{\alpha,\beta}X\longrightarrow
\P_{g(\alpha),g(\beta)}Y$ is a trivial fibration of simplicial
sets. So it satisfies the RLP with respect to any trivial cofibration of
simplicial sets. Therefore $f=\glob(i)$ satisfies the LLP with respect
to $g$ by Lemma~\ref{llp1}. So $\mathcal{J}\subset \cof(\mathcal{I})$. Since
$\cof(\mathcal{I})$ is closed under pushouts and transfinite
compositions, one deduces the inclusion $\cell(\mathcal{J})\subset
\cof(\mathcal{I})$.

The inclusion $\cell(\mathcal{J})\subset \mathcal{W}$ is the consequence
of several facts. Let $f:U\longrightarrow V$ be a morphism of
simplicial sets. Consider the pushout diagram of $\dtop(\sis,\p)$:
\[
\xymatrix{
\glob(U)\fr{} \fd{\glob(f)} & X \fd{g}\\
\glob(V)\fr{} & Y \cocartesien}
\] 
Then by \cite{model3} Proposition~15.2, for any $(\alpha,\beta)\in
X^0\p X^0$, the morphism of simplicial sets
$\P_{\alpha,\beta}X\longrightarrow \P_{g(\alpha),g(\beta)}Y$ is a
transfinite composition of pushouts of morphisms of the form
\[\id\p\dots\p\id\p f \p\id\p\dots\p\id.\] 
If $f$ is a (generating) trivial cofibration of $\sis$, then
$\id\p\dots\p\id\p f \p\id\p\dots\p\id$ is a trivial cofibration as
well since any object of the monoidal model category $(\sis,\p)$ is
cofibrant. So the morphism of simplicial sets
$\P_{\alpha,\beta}X\longrightarrow \P_{g(\alpha),g(\beta)}Y$ is a
trivial cofibration of simplicial sets, and in particular a weak
homotopy equivalence. Therefore $\cell(\mathcal{J})\subset
\mathcal{W}$.

Let $f:X\longrightarrow Y\in\inj(\mathcal{I})$. Then $f^0$ satisfies the
RLP with respect to both $R$ and $C$. So $f^0$ is a bijection of sets
by Lemma~\ref{llp2}. And for any $(\alpha,\beta)\in X^0\p X^0$, the
morphism of simplicial sets $\P_{\alpha,\beta}X\longrightarrow
\P_{f(\alpha),f(\beta)}Y$ satisfies the RLP with respect to any
cofibration of simplicial sets by Lemma~\ref{llp1}. So the morphism of
simplicial sets $\P_{\alpha,\beta}X\longrightarrow
\P_{f(\alpha),f(\beta)}Y$ is a trivial fibration of simplicial
sets. And the morphism of simplicial sets
$\P_{\alpha,\beta}X\longrightarrow \P_{f(\alpha),f(\beta)}Y$ satisfies
the RLP with respect to any trivial cofibration of simplicial
sets. Hence $\inj(\mathcal{I}) \subset
\inj(\mathcal{J})\cap\mathcal{W}$.

Let $f:X\longrightarrow Y\in \inj(\mathcal{J})\cap\mathcal{W}$. Then for
any $(\alpha,\beta)\in X^0\p X^0$, the morphism of simplicial sets
$\P_{\alpha,\beta}X\longrightarrow
\P_{f(\alpha),f(\beta)}Y$ satisfies the RLP with respect to any
trivial cofibration of simplicial sets by Lemma~\ref{llp1}. So the
morphism of simplicial sets $\P_{\alpha,\beta}X\longrightarrow
\P_{f(\alpha),f(\beta)}Y$ is a fibration of simplicial sets, and a
weak homotopy equivalence since $f\in \mathcal{W}$. Therefore for any
$(\alpha,\beta)\in X^0\p X^0$, the morphism of simplicial sets
$\P_{\alpha,\beta}X\longrightarrow
\P_{f(\alpha),f(\beta)}Y$ satisfies the RLP with respect to any 
cofibration of simplicial sets.  Hence the inclusion
$\inj(\mathcal{J})\cap\mathcal{W}\subset\inj(\mathcal{I})$ by
Lemma~\ref{llp1}.

So far, we have proved that $\dtop(\sis,\p)$ is a cofibrantly
generated model category. It remains to prove that the Quillen
equivalence $|-|:\sis\rightleftarrows \top : \sing$ gives rise to a
Quillen equivalence $\dtop(\sis,\p)\rightleftarrows
\dtop(\top,\p)$. Since the geometric realization functor commutes with
binary products, it gives rise to a well-defined functor from
$\dtop(\sis,\p)$ to $\dtop(\top,\p)$. Since the singular nerve functor
is a right adjoint, it commutes with binary products as well. So it
gives rise to a well-defined functor from $\dtop(\top,\p)$ to
$\dtop(\sis,\p)$. It is routine to prove that this pair of functors
defines an adjunction between $\dtop(\sis,\p)$ and $\dtop(\top,\p)$.

A morphism $f:X\longrightarrow Y$ of either $\dtop(\sis,\p)$ or
$\dtop(\top,\p)$ is a fibration (resp. a trivial fibration) if and
only if for any $(\alpha,\beta)\in X^0\p X^0$, the morphism
$\P_{\alpha,\beta}X\longrightarrow
\P_{f(\alpha),f(\beta)}Y$ is a fibration (resp. a trivial fibration and 
$f^0$ is a bijection of sets). So the functor from $\dtop(\top,\p)$ to
$\dtop(\sis,\p)$ preserves fibrations and trivial
fibrations. Therefore it is a right Quillen functor and the adjunction
is actually a Quillen adjunction between $\dtop(\sis,\p)$ and
$\dtop(\top,\p)$.

It remains to prove that this Quillen adjunction is actually a Quillen
equivalence. Let $X$ be a cofibrant object of $\dtop(\sis,\p)$. Let
$Y$ be a fibrant object of $\dtop(\top,\p)$ (that is, any object of
$\dtop(\top,\p)$). One has to prove that $|X|\longrightarrow Y$ is a
weak equivalence of $\dtop(\top,\p)$ if and only if $X\longrightarrow
\sing Y$ is a weak equivalence of $\dtop(\sis,\p)$. Since any weak
equivalence in $\dtop(\sis,\p)$ or $\dtop(\top,\p)$ induces a
bijection between the $0$-skeletons, it remains to prove that for any
$\alpha,\beta\in X^0=Y^0$, the continuous map
$|\P_{\alpha,\beta}X|\longrightarrow \P_{\alpha,\beta}Y$ is a weak
homotopy equivalence of topological spaces if and only if the morphism
of simplicial sets $\P_{\alpha,\beta}X\longrightarrow
\sing\P_{\alpha,\beta}Y$ is a weak homotopy equivalence of simplicial sets. 
Since any simplicial set is cofibrant, and since any topological space
is fibrant, this follows from the fact that the pair of functors
$(|-|,\sing)$ gives rise to a Quillen equivalence between $\sis$ and
$\top$. \epf

Considering $\dtop(\sis,\p)$ does not prove
Theorem~\ref{verscellulaire} and does not allow to take in account the
T-homotopy equivalences.  Indeed, one still has:

\bth Let $({\rm Cof},{\rm Fib},\mathcal{W})$ be a model structure on
$\dtop(\sis,\p)$ such that a morphism of $\mathcal{W}$ is never a
non-trivial pushout of $R:\{0,1\}\longrightarrow \{0\}$. Then: 1) $R$
is necessarily a cofibration, and therefore there exists a cofibration
which is not a monomorphism; 2) any weak equivalence of $\mathcal{W}$
induces a bijection of sets between the $0$-skeletons.
\eth

\bpf 
The proof goes exactly as the proofs of Lemma~\ref{main00} and
Theorem~\ref{main}. \epf

We are going to use a little bit the theory of locally presentable
categories.  A good reference is \cite{MR95j:18001}.

\bp The category $\dtop(\sis,\p)$ is locally finitely presentable. \ep

\bpf Let $\mathbb{S}$ be the limit-sketch corresponding to the category 
of ``small categories without identities''. The category of models
$\Mod(\mathbb{S},\sis)$ of the sketch $\mathbb{S}$ in the category
of simplicial sets is locally finitely presentable by
\cite{MR95j:18001} Theorem~1.53 since the category
of simplicial sets is locally finitely presentable. That does not
complete the proof because the category $\Mod(\mathbb{S},\sis)$
corresponds to flows having a $0$-skeleton which is not necessarily
a discrete simplicial set anymore.

Now consider the adjunction $\pi_0:\sis\leftrightarrows \set:D$
between simplicial sets and sets where $\pi_0$ is the path-connected
component functor and where for any set $S$, $D(S)$ is the discrete
simplicial set associated to $S$. Then $D:\set\longrightarrow \sis$
can be extended in an obvious way to a functor
$\widehat{D}:\dtop(\sis,\p)\longrightarrow
\Mod(\mathbb{S},\sis)$. The functor $\widehat{D}$ is
limit-preserving and colimit-preserving since any limit and colimit of
discrete simplicial sets is a discrete simplicial set.

Let $Z$ be an object of $\Mod(\mathbb{S},\sis)$. Then any morphism
$Z\longrightarrow \widehat{D} X$ factors as a composite
$Z\longrightarrow \widehat{D}T\longrightarrow \widehat{D} X$ where the
cardinal of the underlying set of $T$ is lower than the cardinal of
the underlying set of $Z$. Thus there exists a set of solutions which
proves that $\widehat{D}$ admits a left adjoint
$\widehat{\pi_0}:\Mod(\mathbb{S},\sis)\longrightarrow
\dtop(\sis,\p)$ by Freyd's adjoint functor theorem. 
Since the category $\sis$ (resp. $\set$) can be viewed as a full
subcategory of $\Mod(\mathbb{S},\sis)$ (resp. $\dtop(\sis,\p)$) by the
$0$-skeleton, $\widehat{\pi_0}(K)=\pi_0(K)$ for any simplicial set
$K$. Moreover, one has $\widehat{\pi_0}\circ \widehat{D}=\id$.

It is already known that $\dtop(\sis,\p)$ is cocomplete. Let
$X\in\dtop(\sis,\p)$. Then $\widehat{D}(X)$ is isomorphic to a
directed colimit $\liminj X_i$ where the $X_i$ are finitely
presentable. Then $X\iso \widehat{\pi_0}(\widehat{D}(X))\iso \liminj
\widehat{\pi_0}(X_i)$ since the functor $\widehat{\pi_0}$ is a left
adjoint. It then suffices to prove that the $\widehat{\pi_0}(X_i)$ are
finitely presentable flows.

Let $\liminj Y_j$ be a directed colimit of $\dtop(\sis,\p)$. Then 
\begin{align*}
& \dtop(\sis,\p)(\widehat{\pi_0}(X_i),\liminj Y_j) & \\
& \iso \Mod(\mathbb{S},\sis)(X_i,\widehat{D}(\liminj Y_j)) & \hbox{ by adjunction}\\
& \iso \Mod(\mathbb{S},\sis)(X_i,\liminj \widehat{D}(Y_j)) & \hbox{ since $\widehat{D}$ colimit-preserving}\\
& \iso \liminj \Mod(\mathbb{S},\sis)(X_i,\widehat{D}(Y_j)) & \hbox{ since $X_i$ finitely presentable}\\
& \iso \liminj \dtop(\sis,\p)(\widehat{\pi_0}(X_i),Y_j) & \hbox{ by adjunction again.}
\\
\end{align*}
\epf

\bpf[Proof of Theorem~\ref{verscellulaire}]
By Theorem~\ref{modelsis}, the model category $\dtop(\top,\p)$ is
Quillen equivalent to the locally presentable model category
$\dtop(\sis,\p)$. By Dugger's works \cite{MR2002k:18022} and
\cite{MR2002k:18021}, the model category $\dtop(\sis,\p)$ is Quillen
equivalent to some Bousfield localization of the model category
consisting of the simplicial presheaves over the cofibrant
$\lambda$-presentable objects of $\dtop(\sis,\p)$ for a regular
cardinal $\lambda$ which must be large enough and equipped with the
Bousfield-Kan model structure \cite{MR51:1825}.  The latter turns out
to be cellular. Therefore all its cofibrations are effective
monomorphisms. \epf

To write down the proof of Theorem~\ref{verscellulaire}, we used a
locally presentable model category which is Quillen equivalent to the
model category of topological spaces. Other choices were possible.
Instead of considering the category of simplicial sets, it would have
been possible for instance to take the category of $\Delta$-generated
topological spaces \cite{notedelta}. This category is the largest full
subcategory of $\top$ such that the full subcategory of
$n$-dimensional simplices $\Delta^n$ for $n\geq 0$ is dense in
it. Since the full subcategory of simplices is small, Vop\v{e}nka's
principle ensures that this category is locally presentable. Moreover,
there exists an unpublished proof of this fact due to Jeff Smith which
does not make use of Vop\v{e}nka's principle
\cite{papier}.

The weak S-homotopy model category of flows is replaced by another one
Quillen equivalent and cellular. What do the flows with empty space
become in the new model category ? The restriction of the weak
S-homotopy model category of flows to the category of sets is the
model structure $({\rm Cof},{\rm
Fib},\mathcal{W})=(\All,\All,\Iso)$. It is locally presentable and the
path-connected component functor $\pi_0:\sis\longrightarrow\set$
induces a homotopically surjective morphism of model categories in the
sense of
\cite{MR2002k:18022}. So there exists a set $S$ of morphisms of simplicial 
sets such that Bousfield-localizing $\sis$ by $S$ makes $\pi_0$ a Quillen
equivalence. It is actually possible to take $S=\{\mathbf{S}^1\subset
\mathbf{D}^2\}$ since in the latter model structure, a weak
equivalence between two simplicial sets is indeed a morphism of
simplicial sets inducing a bijection between the path-connected
components (\cite{ref_model2} Section~1.5 ``Postnikov
approximations'').  So the answer is: the sets are replaced by the
simplicial sets and the epimorphism $R:\{0,1\}\longrightarrow \{0\}$
is replaced by the effective monomorphism $\{0,1\}\subset [0,1]$.

\section{Concluding discussion}\label{conc}

The main result of this paper is that the category of flows cannot be
the underlying category of a model category whose corresponding
homotopy types are the flows up to weak dihomotopy. A hint of how the
underlying category of such model category could look like is given by
Dugger's works on combinatorial model categories. The new candidate
for the study of dihomotopy is another model category whose
cofibrations are effective monomorphisms. It contains more objects and
the same weak S-homotopy types. For example, the objects without path
space are the simplicial sets, not only the discrete ones. This way,
the problems existing because of the badly-behaved cofibration
$R:\{0,1\}\longrightarrow\{0\}$ do not appear anymore. It remains to
see whether the Bousfield localization of this new model category with
respect to T-homotopy equivalences has the correct behaviour.

\end{document}